\documentclass[reqno]{amsart}
\usepackage{amsmath,amsbsy,amsfonts}
\usepackage{color}
\date{}
\usepackage[notref,notcite]{showkeys}
\def\nd{\noindent}
\def\thend{\rule{3mm}{3mm}}

\newtheorem{theorem}{Theorem}[section]

\newtheorem{prop}{Proposition}[section]
\newtheorem{lem}{Lemma}[section]
\newtheorem{rmk}{Remark}[section]

\newcommand{\w}{W_0^{1,\Phi}(\Omega)}

\newcommand{\cqd}{\hspace{10pt}\fbox{}}

\newcommand{\Int}{\displaystyle\int_{\Omega}}

\newcommand{\ds}{\displaystyle}
\newcommand{\Fr}{\displaystyle\frac}
\newcommand{\Proof}{{\hspace{-.8cm}\bf Proof: }}

\newcommand{\on}{o_n(1)}
\newcommand{\el}{\ell^*}

\newcommand{\R}{\mathbb{R}}

\newcommand{\N}{\mathcal{N}}

\begin{document}
\title[Multiplicity of solutions for a quasilinear  problem] {Multiplicity of solutions for a nonhomogeneous  quasilinear elliptic problem with critical growth}
\vspace{1cm}

\author{M. L. M. Carvalho}
\address{M. L. M. Carvalho \newline Universidade Federal de Goi\'as, IME, Goi\^ania-GO, Brazil }
\email{\tt marcos$\_$leandro$\_$carvalho@ufg.br}

\author{J. V. Goncalves}
\address{J. V. Goncalves. \newline Universidade Federal de Goi\'as, IME, Goi\^ania-GO, Brazil }
\email{\tt goncalves.jva@gmail.com}

\author{C. Goulart}
\address{C. Goulart \newline Universidade Federal de Goi\'as, Regionsl Jata\'{\i}, Jata\'{\i}, Brazil }
\email{\tt claudiney@ufg.br}

\author{O. H. Miyagaki}
\address{O. H. Miyagaki \newline Universidade Federal de Juíz de Fora, Juiz de Fora, Brazil }
\email{\tt ohmiyagaki@gmail.com }

\subjclass{35J20, 35J25, 35J60, 35J92, 58E05} \keywords{Variational methods, Quasilinear Elliptic Problems, Nehari method, Sign-changing solutions}
\thanks{O. H. Miyagaki is corresponding author and he received research grants from CNPq/Brazil  304015/2014-8 and INCTMAT/CNPQ/Brazil.}

\begin{abstract}
It is established some existence and multiplicity of solution  results for a quasilinear elliptic problem driven by  $\Phi$-Laplacian operator. One of these solutions is  built as a ground state solution. In order to prove our main results we apply the Nehari method combined with the concentration compactness theorem  in an  Orlicz-Sobolev framework. One of the difficulties in dealing with this kind of operator is the lost of homogeneity properties. 
\end{abstract}

\maketitle

\section{Introduction}

In this work we  will    establish some existence and multiplicity   results for the following quasilinear elliptic problem 
\begin{equation}\label{eq1}
\left\{\begin{array}{rclcl}
-\Delta_{\Phi} u  & = & |u|^{\ell^*-2}u+f,  &\mbox{ in }&    \Omega, \\
u  &=& 0,&  \mbox{ on } &\partial\Omega ,
\end{array} \right.
\end{equation}
where  $\Delta_{\Phi}$ denotes the  $\Phi-$laplacian operator, which is defined by $\Delta_{\Phi} u = \mbox{div} (\phi(|\nabla u|)\nabla u),$ $\ell^* = \ell N/(N - \ell)~ (1< \ell < N),$ $\Omega\subset\mathbb{R}^{N}$ is bounded and smooth domain,  $f \underset{\neq 0}{\geq 0},$  and in order  to simplify the technicalities we assume $f\in L^{\frac{\ell N}{N(\ell-1)+\ell}}(\Omega)=L^{^{\frac{\el}{\el-1}}}(\Omega)\equiv L^{{\el}^{'}}(\Omega).$ With respect to the function $\phi: (0,\infty)\rightarrow (0,\infty),$ we assume that it is $C^{2}$ and satisfies the following conditions
\begin{description}
	\item[$(\phi_1)$] $\displaystyle \lim_{t \rightarrow 0} t \phi(t)= 0, \displaystyle \lim_{t \rightarrow \infty} t \phi(t)= \infty$;
	\item[$(\phi_2)$] $t\mapsto t\phi(t)$ is strictly increasing;
	\item[$(\phi_3)$] $ -1<\ell-2:=\ds\inf_{t>0}\Fr{(t\phi(t))''t}{(t\phi(t))'}\leq \ds\sup_{t>0}\Fr{(t\phi(t))''t}{(t\phi(t))'}=:m-2<N-2.$
\end{description}

Furthermore, we  shall assume the following hypothesis
\begin{flushleft}
$(H)$ $\hspace*{.7cm}$ $1<\Fr{\ell(\el-m)}{\el-\ell}\leq \ell\leq m<\el.$
\end{flushleft}
\begin{rmk}
	Notice that the above inequalities still hold when:
	\begin{enumerate}
		\item $\Phi(t)=pt^{p-2}$ with $1 <p<\infty $ and  $\ell=m=p,$ in this case  $\Delta_\Phi =\Delta_p$, where $\Delta_p$ denotes the  $p-$Laplacian operator.
		\item $\Phi(t)=pt^{p-2}+ qt^{q-2}$ with $1<p< q <\infty,$ $ \ell=p$ and $m=q,$in this case $\Delta_\Phi$ turns  the $\Delta_p +\Delta_q$ operator.  Here $\Delta_p +\Delta_q$   denotes the  $(p,q)-$ laplacian operator. See 	\cite{MugnaiPapageorgiou, tanaka})  for this kind of operators.	
		\item Other examples, for instance involving anisotropic elliptic problems, can be seen in \cite{CSG} and references therein.
	\end{enumerate}

\end{rmk}

The main difficulty in dealing with this kind of operator is because it is inhomogeneous, which requires some aditional effort to overcome the estimates. As is mentioned in \cite{fang} the problem has many physical applications, for instance, in  nonlinear elasticity, plasticity, generalized Newtonian fluids, etc. 
We refer the reader to the following   related papers \cite{Bouchekif,giovanyquoirin, Fuk_1,Fuk_2,Motreanu, fang} and in references therein, where there have handled handled  different types of  nonlinearities involving this kind of  operator. Problems like above was started in a beautiful work due to Br\'ezis and Nirenberg\cite{BN2},  when $\Delta_\Phi=\Delta$, where they treated a nonhomogeneous problem with critical growth obtaining existence result, assuming that   $f \underset{\neq 0}{\geq 0},$ together with some aditional conditions.  Then Tarantello \cite{tarantello}  treated the same problem getting  existence and multiplicity results   under a stronger hypothesis that made in \cite{BN2}. These works were extended in  \cite{hirano},which was obtained four weak solutions, at least one of them is sign
changing solution. 	On the other hand, in \cite{clapp} is proved some multiplicity results for symmetric domain by using the category theory. There are only few works involving $p-$ Laplacian, that is, when  $\Delta_\Phi=\Delta_p,$ extending results in \cite{tarantello}. We would like to mention \cite{Cao,Chabro} and  references therein.

Due to the nature of the operator
$\displaystyle \Delta_{\Phi} $ we shall work in the framework of Orlicz-Sobolev spaces $W^{1,\Phi}_{0}(\Omega)$. Throughout this paper we define
$$
\Phi(t) = \int_{0}^{t} s \phi(s) ds, t \geq 0,
$$
which is extended as even function, $\Phi(t)=\Phi(-t),$ for all $t<0.$

Recall that hypotheses $(\phi_1)-(\phi_2)$ allow us to use the Orlicz and Orlicz-Sobolev spaces, while the hypothesis $(\phi_3)$ ensures that
the Orlicz-Sobolev spaces are Banach reflexive spaces. There are several publications on  Orlicz-Sobolev spaces, we would like to recommend  the  reader  to
\cite{A,DT, Fuk_1, gossez-Czech, Rao1,fang}. However, for the  sake of  completeness, we recall some definitions and properties in the Appendix.

\vskip.2cm

From the continuous embedding  $W_0^{1,\Phi}(\Omega)\hookrightarrow L^{\el}(\Omega)$,(see \cite{A,DT}), we define 
\begin{equation}\label{melhor}
S=\displaystyle\inf\left\{\Fr{||u||^\ell}{||u||^\ell_{\el}},\ \ u\in W^{1,\phi}(\Omega)\setminus\{0\}\right\}.
\end{equation}

Since our approach is variational method, the  functional  $J : W^{1, \Phi}_{0}(\Omega) \rightarrow \mathbb{R}$ associated with our problem is given by
$$
J(u)=\Int\Phi(|\nabla u|)-\Fr{1}{\el} \Int |u|^{\el}- \Int f u,\ \ u\in W_0^{1,\Phi}(\Omega),
$$
is well-defined and of class $C^{1}$. The Euler-Lagrange equations for $J$ are precisely the weak solutions for problem \eqref{eq1}. Hence finding weak solutions for the problem \eqref{eq1} is equivalent to find critical points for the functional $J$. Here we emphasize that $J$ is in $C^{1}$ class due the hypotheses built on the function $\phi$. This is the main reason in order to consider the hypothesis $(\phi_{3})$ that is crucial in our arguments. The Gateaux derivative for $J$ possesses the following form
\begin{equation}\label{functional}
\left<J'(u),v\right> = \Int\phi(|\nabla u|)\nabla u\nabla v- \Int |u|^{\el-2}uv- \Int f v\nonumber
\end{equation}
for any $u,v\in W_0^{1,\Phi}(\Omega)$. In general, using hypotheses $(\phi_{1})-(\phi_{3})$, the functional $J$ is not in $C^{2}$ class.

In order to perfom our precise hypotheses for our results, we will consider the functions  $	g_{\alpha}:[0,\infty)\rightarrow \mathbb{R},~\alpha\in\{\ell,m\}$ defined by 
\begin{equation}\label{galpha}
	g_{\alpha}(t):=g_{\alpha,u}(t)=t^{\alpha-1}\Int \phi(|\nabla u|)|\nabla u|^2-t^{\el-1}|u|^{\el},\ \ t>0.
\end{equation}
It is easy to see that there exists $\overline{t}_\alpha>0$ such that
  $$g_{\alpha}(\overline{t}_\alpha)=\displaystyle \max_{t>0}g_{\alpha}(t).$$

Inspired by \cite{tarantellonewmann}, given $u\in W^{1,\phi}_{0}(\Omega),$ with $ \ \ ||u||_{{\el}}=1,$ we assume the following assumptions on $ f. $

\begin{description}
\item[$(f_1)$] Suppose either $\overline{t}_\ell,\overline{t}_m\geq 1$ or $\overline{t}_\ell,\overline{t}_m\leq 1.$ Then 
$$||f||_{({\el})'}\leq \lambda_1:= \min\left\{
S^{\frac{\alpha(\el-1)}{\ell(\el-\alpha)}}\left[\Fr{\ell(\ell-1)}{\el-1}\right]^{\frac{\alpha-1}{\el-\alpha}}
	\left[\Fr{\ell(\el-m)}{\el-1}\right]: ~\alpha=\ell,m\right\};$$
	\item[$(f_2)$] If $\ell<m$ and $\overline{t}_\ell\leq 1\leq \overline {t}_{m}$ hold, we suppose  $$||f||_{({\el})'}\leq \min\left\{\lambda_1,\Fr{\el-m}{m-1}\right\}.\mbox{(See Lemma \ref{nehari-})}$$
\end{description}
 We have a second solution to the problem \eqref{eq1} considering a more restrictive condition given by:
\begin{description}
\item[$(f_2)'$] If $\ell<m$ and $\overline{t}_\ell\leq 1\leq \overline {t}_{m}$ hold, we assume $$||f||_{({\el})'}\leq \min\left\{\frac{1}{m}\lambda_1,\Fr{\el-m}{m-1}\right\}.(\mbox{
	See Lemma \ref{nehari-}})$$\end{description}

Our first main result can be read as follows
\begin{theorem}\label{teorem1}
In addition to	 $(\phi_{1}) - (\phi_{3})$ and $(H)$, suppose  $f \underset{\neq 0}{\geq 0},$  and  $f\in L^{{\el}^{'}}(\Omega).$ Assume either $(f_1)$ or $(f_2)$ holds. Then there exists $\Lambda_1> 0$ such that problem \eqref{eq1} admits at least one positive ground state solution $u^+$ satisfying $J(u^{+}) \leq 0$ for any $f$ such that $0 < ||f||_{(\el)'} < \Lambda_{1}$.
\end{theorem}

Now we shall consider the following result

\begin{theorem} \label{teorema2}
Suppose $(\phi_{1}) - (\phi_{3})$ and  $(H).$ Assume  $f \underset{\neq 0}{\geq 0},$  and  $f\in L^{{\el}^{'}}(\Omega),$ and  either $(f_1)$ or $(f_2)'$  holds. Then there exists $\Lambda_2 > 0$ in such way that problem \eqref{eq1} admits at least one positive  solution $u^{-}$ satisfying $J(u^{-}) > 0$ for any  $f$ verifying $0 <||f||_{(\el)'} < \Lambda_{2}$.
\end{theorem}

Putting together the all results established just above and using a regularity result for quasilinear elliptic problems we can state the following multiplicity result.

\begin{theorem}\label{teorema3}
In addition to  $(\phi_{1}) - (\phi_{3})$ and $(H)$, suppose  $f \underset{\neq 0}{\geq 0},$  and  $f\in L^{{\el}^{'}}(\Omega).$  Assume either $(f_1)$ or $(f_2)'$  holds. Then problem \eqref{eq1} admits at least two positive $u^{+}, u^{-}$ which belong to $C^{1,\alpha}(\bar\Omega)$ whenever $0 < ||f||_{(\el)'} < \Lambda=\min\{\Lambda_1,\Lambda_2\}$. Furthermore, the function
$u^+$ is a ground state solution for each $f$ satisfying $0 < ||f||_{(\el)'} < \Lambda$.
\end{theorem}

\begin{rmk}
We point out that concerning just existence of solution,  $f$ can change sign, see Lemma \ref{fib}. However in such case the solution could change sign, as well.
\end{rmk}

\section{Preliminary results}

In this section we give some basic results involving  the Nehari manifold method, including the fibering maps  associated with  the  functional $J,$ which will give information on the critical points of  Euler-Lagrange functional $J$.  We suggest the reader to the book due to Willem \cite{Willem}, for an overview on the Nehari method. The proofs of our results follow closely the  arguments used in \cite{CSGG,JME}.

The Nehari manifold associated with the functional $J$ is
given by
\begin{equation}
  \begin{array}{rcl}\label{nehari}
	\mathcal{N}&=&\{u\in W_0^{1,\Phi}(\Omega)\setminus \{0\}: \left<J'(u),u\right>=0\}\\[2ex]
  \end{array}
\end{equation}
It will be proved later on  that $\mathcal{N}$ is a $\mathcal{C}^1$-submanifold of $W^{1,\Phi}_{0}(\Omega)$.

Initially, note that if $u\in \mathcal{N},$ by \eqref{nehari}, we have that

\begin{equation}\label{eq2}
J(u)=\Int \Phi(|\nabla u|) -  \phi(|\nabla u|)|\nabla u|^2+\left(1-\Fr{1}{\el}\right)|u|^{\el},
\end{equation}
or equivalently
\begin{equation}\label{eq3}
J(u)=\Int \Phi(|\nabla u|) - \Fr{1}{\el} \phi(|\nabla u|)|\nabla u|^2-\left(1-\Fr{1}{\el}\right)fu.
\end{equation}

First of all  we shall prove some geometric properties of functional  $J,$ which allows us to find  a critical point for $J$. 
\begin{prop}\label{coercive}
The functional $J$ is coercive and bounded from  below on $\mathcal{N}$.
\end{prop}
\proof In virtue of $(\phi_3)$, we have  $m\Phi(t)\geq t^2\phi(t)$ for each $t \geq 0$. Using this fact and \eqref{eq3}, we obtain
\begin{equation}\label{eq4}
J(u)\geq\left(\Fr{1}{m}-\Fr{1}{\el} \right)\Int\phi(|\nabla u|)|\nabla u|^2+\left(\Fr{1}{\el}-1\right)\Int fu.\nonumber
\end{equation}

Now by combining 
$$\min\{||u||^{\ell},||u||^m\}\leq\Int \Phi(|\nabla u|)\leq \Fr{1}{\ell}\Int\phi(|\nabla u|)|\nabla u|^2$$
with the H\"older inequality and the continuous embedding   $W_0^{1,\Phi}(\Omega)\hookrightarrow L_{\Phi_*}(\Omega)\hookrightarrow L^{\ell^*}(\Omega)$, we obtain
\begin{equation}\label{functional1}
J(u)\geq \ell\left(\Fr{1}{m}-\Fr{1}{\el} \right)\min\{||u||^\ell,||u||^m\}+S^{\frac{-1}{\ell}}\left(\Fr{1}{\el}-1\right)||f||_{(\el)'}||u||,
\end{equation} where $S$ is given by \eqref{melhor}. 
Thus, $J$ is coercive and bounded  from  below on $\mathcal{N}$. The proposition is proved.\hfill\cqd

Now, define the fibering map $\gamma_u: (0,+\infty)\to \mathbb{R}$ given by 
$$\gamma_u(t):=J(tu)=\Int\Phi(t|\nabla u|)-\Fr{t^{\el}}{\el}|u|^{\el}- tfu.$$
From $(\phi_{1}) - (\phi_{2})$ it follows that $\gamma_{u}$ is of $C^{1},$ and  its Gateaux derivative is given by

\begin{equation}\label{priderivada}\gamma'_u(t)= t\Int \phi(t|\nabla u|)|\nabla u|^2 -t^{\el-1} |u|^{\el}-fu.\end{equation}
 The main feature of the fibering map is the knowledge of the geometry of $\gamma_{u}$, which will  give information about the existence and multiplicity of solutions. This method was introduced in \cite{Drabek}, then it was also employed, for instance,  in \cite{Brown2,Brown1,Brown3,tarantello,tarantellonewmann,wu,wu3} and references therein.

\begin{rmk}
Notice that $tu\in\mathcal{N}$ if, and  only if, $\gamma'_u(t)=0.$ Therefore, $u\in\mathcal{N}$ if, and  only if, $\gamma'_u(1)=0.$ Thus, the  stationary points of fibering map are the critical points of $J$ on $\mathcal{N}$.
\end{rmk}

Define $\psi(u)=\left<J'(u),u\right>,\  u \in W_0^{1,\Phi}(\Omega).$ Then,  for all  $u  \in W^{1,\Phi}_{0}(\Omega),$ we have

\begin{equation}\label{psi}
\left<\psi'(u),u\right>=\Int \phi'(|\nabla u|)|\nabla u|^3+2\phi(|\nabla u|)|\nabla u|^2-\el |u|^{\el}- fu.
\end{equation}

As was made in  Tarantello in \cite{tarantello,tarantellonewmann}, let us  split $\N$ into three sets, namely,
$$\N^+:=\{u\in \N:\left<\psi'(u),u\right>>0\};$$
$$\N^-:=\{u\in \N:\left<\psi'(u),u\right><0\};$$
$$\N^0:=\{u\in \N:\left<\psi'(u),u\right>=0\},$$
  which correspond to the  critical points of minimum, maximum and inflexions points, respectively.

\begin{rmk}\label{gamma''}
For $u\in \N,$ by \eqref{eq2} and \eqref{eq3}, we have \begin{equation}\label{eq5}\begin{array}{rcl}\left<\psi'(u),u\right>&=& \Int [\phi'(|\nabla u|)|\nabla u|^3+\phi(|\nabla u|)|\nabla u|^2]-(\el-1)|u|^{\el}\\[2ex]
&=&
\Int [\phi'(|\nabla u|)|\nabla u|^3+(2-\el)\phi(|\nabla u|)|\nabla u|^2]-(1-\el)fu.
\end{array}\end{equation}
\end{rmk}

The next result is the crucial step in our argument to prove the main result.

\begin{lem}\label{c1} 
Suppose either $(f_1)$ or $(f_2),$ and $(\phi_1)$-$(\phi_3)$ hold. Then,
  
\begin{enumerate}
\item $\N^0=\emptyset$.
\item $\N=\N^+\dot{\cup} \ \N^-$ is a $C^1$-manifold.
\end{enumerate}\end{lem}
\proof {\bf Proof of item (1).} Assume by contradiction  that $\N^0\neq\emptyset.$ Fix $u\in\N^0.$  Then,  $\gamma'_u(1)=\left<\psi'(u),u\right>=0.$ From \eqref{nehari} and \eqref{eq5}, we obtain, $$0=\left<\psi'(u),u\right>=\Int \phi(|\nabla u| )|\nabla u|^2+\phi'(|\nabla u|)|\nabla u|^3+(1-\el)|u|^{\el}.$$ By hypothesis $(\phi_3)$ we infer that
$$(\ell-1)\Int \phi(|\nabla u| )|\nabla u|^2\leq(\el-1)||u||_{\el}^{\el}
\leq (\el-1)S^{-\frac{\el}{\ell}}||u||^{\el},$$ where $S$ is the best constant of the embedding $W_0^{1,\Phi}(\Omega)\hookrightarrow L^{\el}(\Omega).$ 
 On the other hand,
$$(\ell-1)\Int \phi(|\nabla u| )|\nabla u|^2\geq\ell(\ell-1)\Int \Phi(|\nabla u|)\geq \ell(\ell-1)\min\{||u||^\ell,||u||^m\}. $$
Comparing the above two expressions, we conclude that 

\begin{equation}\label{des1}
	||u||\geq \left[\Fr{\ell(\ell-1)S^{\frac{\ell^*}{\ell}}}{(\el-1)}\right]^{\frac{1}{\el-\alpha}}.
\end{equation}
Now, using \eqref{eq5}, we get
$$0=\left<\psi'(u),u\right>=\Int (2-\el)\phi(|\nabla u| )|\nabla u|^2+\phi'(|\nabla u|)||\nabla u|^3+(\el-1)fu.$$ 
From $(\phi_3)$, we obtain
$$(\el-m)\Int \phi(|\nabla u| )|\nabla u|^2\leq(\el-1)\Int fu$$
 Arguing as above, we get
$$\ell(\el-m)\min\{||u||^\ell,||u||^m\}\leq (\el-1)\Int fu.$$
Therefore, from the H\"older's inequality, we get
\begin{equation}\label{des2}\min\{||u||^\ell,||u||^m\}\leq \Fr{(\el-1)}{\ell(\el-m)}\Int fu \leq \Fr{S^{-\frac{1}{\ell}}(\el-1)}{\ell(\el-m)} ||f||_{(\el)'}||u||.\end{equation}
Comparing \eqref{des1} and \eqref{des2}, we get
\begin{equation}\label{lambda1}
||f||_{(\el)'}\geq S^{\frac{\alpha(\el-1)}{\ell(\el-\alpha)}}\left[\Fr{\ell(\ell-1)}{\el-1}\right]^{\frac{\alpha-1}{\el-\alpha}}
	\left[\Fr{\ell(\el-m)}{\el-1}\right]\geq\lambda_1, ~\alpha=\ell,m,
\end{equation}
which is a contradiction if we assume either $(f_1)$ or $(f_2)$.
  \vspace{.2cm}

\noindent {\bf Proof of item (2).} Suppose without loss of generality that, $u\in \N^+$. 

Define $G(u):=\left<J'(u),u\right>.$ We can see that
 $$G'(u)=\left<J''(u)\cdot(u,u)\right>+\left<J'(u),u\right>=\left<\psi'(u),u\right> >0,\,\,\ \forall u\in\N^+.$$  
Furthermore, using  \eqref{nehari}, we also have that $\left<J'(u),u\right>=0$. Hence, $0\in\R$  is a regular value for $G$ and $\N^+=G^{-1}(0).$  That is, $\N^+$ is a $C^1$-manifold. Similarly , we may show that $\N^-$ is a $C^1$-manifold. Hence, since we are supposing  $(f_1)$ and $(f_2),$ the proof of item (2) follows in virtue of  $\N^{0}=\emptyset.$ 

\hfill\cqd
\vskip.2cm

Next we are going to prove that any critical point for $J$ on $\mathcal{N}_\lambda$ is a free critical point, i.e, is a critical point in the whole space $W^{1,\Phi}_{0}(\Omega)$. Actually, the proof of the Lemma below is fairly standard and we include it for the sake of completeness.
\begin{lem}\label{criticalpoint}
Let $u_0$  be a local minimum (or local maximum) of $J.$ If $u_0 \notin$ $\N^0,$ then $u_0$ is a critical point of $J.$.
\proof Suppose without any loss of generality that $u_0$ is a local minimum of $J.$ 
Define the function $$\theta(u)=\left<J'(u),u\right>=\Int\phi(|\nabla u|)|\nabla u|^2-|u|^{\el}-fu.$$ 

Then $u_0$ is a solution for the minimization problem
\begin{equation}\label{lagrange} \min\left\{
 J(u),\  \ 
\theta(u)=0
\right\}.\end{equation}
Proceeding as in Carvalho et al. \cite{JME}, we have 
$$\left<\theta'(u),v\right>=\Int \phi'(|\nabla u|)|\nabla u|\nabla u\nabla v+2\phi(|\nabla u|)\nabla u\nabla v-f v-\el |u|^{\el-2}uv$$ holds true for all  $u, v  \in W^{1,\Phi}_{0}(\Omega)$. Making $u = v = u_{0},$  since $u_0\in \N^{+}$, by \eqref{nehari} and \eqref{eq5}, we get
$$\begin{array}{rcl}
\left<\theta'(u_0),u_0\right>&=&\Int \phi'(|\nabla u_0|)|\nabla u_0|^3+\phi(|\nabla u_0|)|\nabla u_0|^2\\[2ex]
&-&(\el-1)\Int {|u_0|}^{\el}=\left<\psi'(u_0),u_0\right>.\end{array}$$

From Lemma \ref{c1}, the problem \eqref{lagrange}  has a solution  verifying
 $$\left<J'(u_0),u_0\right>=\mu\left<\theta'(u_0),u_0\right>=0,$$
where $\mu\in\R$ which is given by Lagrange multipliers Theorem. Notice that   $\left<\theta'(u_0),u_0\right>\neq 0,$ then $\mu=0$, i.e, $u_{0}$
is a critical point for $J$ on $W^{1,\Phi}_{0}(\Omega)$. The proof  of  lemma is  complete.\hfill\cqd
\end{lem}


Now we give a complete description on the geometry for the fibering map associated
with  problem \eqref{eq1}, where we will foccus on the sign of $\Int f u .$

  Consider the auxiliary function
$$m_u(t)=\Int t\phi(t|\nabla u|)|\nabla u|^2-t^{\el-1}|u|^{\el},$$

\noindent  where the points $t u \in \mathcal{N}$ will compared with  the function $m_{u}$. 

\begin{lem}\label{m_uegamma_u}
Let $t>0$ be fixed. Then $tu\in \N$ if, and only if, $t$ is a solution of $m_u(t)=\Int  fu.$
\proof Fix $t>0$ in such may that $tu\in\N$. Then

$$t\Int\phi(|\nabla t u|)|\nabla u|^2 -t^{\el-1}\Int b(x)|u|^{\el}= \Int f u.$$
From the  definition of $m_u$, the proof of the result follows.\hfill\cqd
\end{lem}

The next lemma will give  a precise information on the  function $m_{u}$ and the fibering map.

\begin{lem}\label{m_u-comp}
 There exists an unique critical point for $m_{u}$, i.e, there is an unique point $\tilde{t}>0$ in such way that $m'_u(\tilde{t})=0$.
 Furthermore, we know that $\tilde{t} > 0$ is a global maximum point for $m_{u}$ and $m_{u}(\infty) = - \infty$.
\proof
 Notice that
 $$m'_u(t)=\Int \phi(t|\nabla u|)|\nabla u|^2+t\Int \phi'(|\nabla (tu)|)|\nabla u|^3-(\el-1)t^{\el-2}|u|^{\el}.$$
 Taking into account $(\phi_{3})$ it is easy to verify that
 \begin{equation}\label{ee2}
 \ell-2\leq\Fr{\phi'(t)t}{\phi(t)}\leq m-2, \,\,\mbox{for any} \,\,\,t \geq 0. \end{equation}

Firstly, we prove that $m_u$ is increasing for $ t > 0$ small enough and $\ds\lim_{t\to\infty}m_u(t)=-\infty$. For $0<t<1,$ using  \eqref{ee2} we get

$$\begin{array}{rcl}
m'_u(t)
&\geq & (\ell-1)t^{m-2}\Int \phi(|\nabla u|)|\nabla u|^2-(\el-1)t^{\el-2}|u|^{\el}\\[2ex]
\end{array}$$ 
Since $m<\el$ we mention that $m'_u(t)>0$ for any $t>0$ small enough. Arguing as above we obtain
$$m_u(t)\leq t^{m-1}\Int \phi(|\nabla u|)|\nabla u|^2-t^{\el-1}|u|^{\el}.$$ Therefore, since $m<\el$, we infer that  $\ds\lim_{t\to\infty}m_u(t)=-\infty.$

Next, we will show that $m_u$ has an unique critical point $\tilde{t}>0.$
Observe that $m'_u(t)=0$ if, and only if, $$t^{2-\el}\Int \phi(t|\nabla u|)|\nabla u|^2+t^{3-\el}\Int \phi'(|\nabla (tu)|)|\nabla u|^3=(\el-1)\Int |u|^{\el}.$$
Define the auxiliary function $\eta_{u}: \mathbb{R} \rightarrow \mathbb{R}$  by
$$\eta_u(t)=t^{2-\el}\Int \phi(t|\nabla u|)|\nabla u|^2+t^{3-\el}\Int \phi'(|\nabla (tu)|)|\nabla u|^3.$$
 Using the inequality below
 $$\eta_u(t)\geq  (\ell-1)t^{m-\el} \Int \phi(|\nabla (u)|)|\nabla u|^2,$$
 it is easy to see that
\begin{equation}\label{ee3}
\ds\lim_{t\to 0^+}\eta_u(t)=+\infty.
\end{equation} 
 
On the other hand, from Proposition \ref{fang}, for any $t>1$, we have

\begin{equation}\label{ee4} 
\eta_u(t)\leq (m-1)t^{m-\el}\Int \phi(|\nabla u|)|\nabla u|^2.
\end{equation}
and
\begin{equation}\label{ee5}
\eta_u(t)\geq (\ell-1)t^{\ell-\el}\Int \phi(|\nabla (u)|)|\nabla u|^2.
\end{equation}
Hence \eqref{ee4} and \eqref{ee5} say that
\begin{equation}\label{ee33}
\lim_{t\to \infty}\eta_u(t)=0.
\end{equation}

\noindent  holds true.

 Moreover, we have also that
$$\begin{array}{rcl}
\eta_u'(t)&=&\Int[(2-\el)t^{1-\el}\phi(t|\nabla u|)|\nabla u|^2+(4-\el)t^{2-\el}\phi'(t|\nabla u|)|\nabla u|^3]+\\[2ex]
&+&t^{3-\el}\Int \phi''(t|\nabla u|)|\nabla u|^4.
\end{array}$$ 
Using hypothesis $(\phi_3)$ we have 
 $$\left\{\begin{array}{rcl}
 t^2\phi''(t)&\leq& (m-4)t\phi'(t)+(m-2)\phi(t),\\
 t^2\phi''(t)&\geq& (\ell-4)t\phi'(t)+(\ell-2)\phi(t).
 \end{array}\right.,$$
 which imply that
 $$\begin{array}{rcl}\eta_u'(t)
&\leq& \Int(m-\el)t^{1-\el}\phi(t|\nabla u|)|\nabla u|^2+\Int (m-\el)t^{2-\el}\phi'(t|\nabla u|)|\nabla u|^3\\
&\leq&(m-\el)(\ell-1) t^{1-\el}\Int\phi(t|\nabla u|)|\nabla u|^2
<0.\end{array}$$
 The proof of this lemma is now complete.
\hfill\cqd
\end{lem}


Next we will estimate $\displaystyle \max_{t>0}m_u(t)$. To do this,  consider $g_\alpha,~\alpha=\ell,m,$ defined in \eqref{galpha}. 
As in the proof of the previous Lemma, there exists  $\overline{t}_{\alpha}>0$, given by 
\begin{equation}\label{tmax}
\overline{t}_\alpha
=\left[\Fr{(\alpha-1)\Int \phi(|\nabla u|)||\nabla u|^2}{(\el-1)\Int |u|^{\el}}\right]^{\frac{1}{\el-\alpha}}>0
\end{equation}
such that  $g_\alpha(\overline{t}_{\alpha})=\displaystyle\max_{t>0}g_\alpha(t).$ 
\begin{rmk}
Notice that  $g_m(t)=g_{\ell}(t)=m_u(t)$ if, only if, $t=1$.
\end{rmk}
\begin{lem}\label{mapf1f2}
	Suppose either $(f_1)$ or $(f_2)$. Then 
	\begin{eqnarray}\label{des-f}
		\max_{t>0}m_u(t)\geq \int_{\Omega}fu,~u\in\w.
	\end{eqnarray}
\end{lem}
\proof If $\displaystyle\int_\Omega fu dx \leq 0$, since $\displaystyle\max_{t>0}m_u(t)>0$, the inequality  \eqref{des-f} is trivially satisfied. Thus, we treat the case   $\displaystyle\int_\Omega fu dx > 0$. Without loss of generality, take  $\|u\|_{\el}=1$ and denote by $A=\Int\phi(|\nabla u|)|\nabla u|^2dx$. 

We will consider three possibilities, namely:
\begin{enumerate}
\item[(i)]: $\overline{t}_{\ell}, \ \overline{t}_{m}\geq 1$. Since
 $\overline{t}_{\ell}\geq 1$, we obtain
\begin{equation}\label{A1}
\el-1\leq(\ell-1)A.
\end{equation}
So that, 
\begin{eqnarray}
1&\leq&\Fr{\ell-1}{\el-1}A 
\leq m\Fr{\ell-1}{\el-1}\max\{||u||^\ell,||u||^m\}\nonumber\\
&=&\max\left\{\left\|\left(m\Fr{\ell-1}{\el-1}\right)^\frac{1}{\ell}u\right\|^\ell,\left\|\left(m\Fr{\ell-1}{\el-1}\right)^\frac{1}{m}u\right\|^m\right\}
= m\frac{\ell-1}{\el-1}\|u\|^m.\nonumber
\end{eqnarray}

On the other hand, using  Proposition \ref{lema_naru} and inequality $$\left(m\Fr{\ell-1}{\el-1}\right)^\frac{1}{\alpha}||u||\geq 1,~\alpha=\ell,m,$$ we get
\begin{eqnarray}
A&\geq& \ell\min\{||u||^\ell,||u||^m\}\nonumber\\
&=&\Fr{\ell(\el-1)}{m(\ell-1)}\min\left\{\left\|\left(m\Fr{\ell-1}{\el-1}\right)^\frac{1}{\ell}u\right\|^\ell,\left\|\left(m\Fr{\ell-1}{\el-1}\right)^\frac{1}{m}u\right\|^m\right\}
=\ell||u||^\ell\nonumber
\end{eqnarray}
\nd Moreover,
\begin{equation*}\label{A2}
\begin{array}{rcl}
\displaystyle\max_{t>0} m_u(t)&\geq& \max g_{\ell}(t)=g_{\ell}(\overline{t}_{\ell})\\[2ex]
&\geq&||u||||u||^{\frac{\el(\ell-1)}{\el-\ell}}\left(\Fr{\ell(\ell-1)}{\el-1}\right)^{\frac{\ell-1}{\el-\ell}}\left(\Fr{\ell(\el-\ell)}{\el-1}\right)\\
&\geq& S^{\frac{(\el-1)\ell}{\ell(\el-\ell)}}\left(\Fr{\ell(\ell-1)}{\el-1}\right)^{\frac{\ell-1}{\el-\ell}}\left(\Fr{\ell(\el-m)}{\el-1}\right)\Fr{1}{||f||_{(\el)'}}\Int fu\\
&\stackrel{(f_1)}\geq& \Int fu.
\end{array}
\end{equation*}

\item[(ii)]: If $\ell<m$ and $\overline{t}_{\ell}\leq 1 \leq \overline{t}_{m}$, then
\begin{equation}\label{B1-1}
\frac{\el-1}{\ell-1}\geq A\geq \frac{\el-1}{m-1}>1.
\end{equation}
Therefore, it follows from  \eqref{B1-1} that
\begin{eqnarray}\label{B2}
\max_{t>0} m_u(t)&\geq& m_u(1) \geq  \frac{\el-m}{m-1}=\frac{\el-m}{m-1}\|u\|_{\el}\nonumber\\
&\geq&\frac{\el-m}{m-1}\frac{1}{\|f\|_{(\el)'}}\Int fu\stackrel{(f_2)}\geq\int_\Omega fu.
\end{eqnarray}
\item[(iii)]: If $\overline t_{\ell}, \ \overline t_{m}\leq 1$, then 
$$(m-1)A\leq\el-1.$$
As in  item  (ii) we get
\begin{equation*}\label{C2}
\begin{array}{rcl}
	\displaystyle\max_{t>0} m_u(t)&\geq& \max g_{m}(t)=g_{m}(\overline t_{m})\\[2ex]
	&\geq& S^{\frac{m(\el-1)}{\ell(\el-m)}}\left(\Fr{\ell(m-1)}{\el-1}\right)^{\frac{m-1}{\el-m}}\left(\Fr{\ell(\el-m)}{\el-1}\right)\Fr{1}{||f||_{(\el)'}}\Int fu\\[2ex]
	&\stackrel{(f_1)}\geq& \Int fu.
\end{array}
\end{equation*}
\end{enumerate}
This finishes the proof of lemma. \hfill\cqd

\begin{lem}\label{fib}
	Let $u\in W^{1,\Phi}_0(\Omega)/\{0\}$ be a fixed function. Then we shall consider the following assertions:
\begin{enumerate}
\item there exists an unique $t_1=t_1(u)> \tilde t$ such that
$\gamma'_u(t_1)= 0$ and $t_1 u\in \N^-$ whenever $\Int fu\leq 0$,
\item suppose either $(f_1)$ or $(f_2)$. Then, if $\Int fu > 0,$ 
there exists unique $0<t_1=t_1(u)<\tilde t<t_2=t_2(u)$ such that 	$\gamma'_u(t_1)=\gamma'_u(t_2)=0$, $ t_1 u\in \N^+$ and $ t_2 u\in \N^-$ .
\end{enumerate}	
\end{lem}
\proof
First of all, notice that  arguing as in \cite{Brown1}, it is easy to see that if  $ tu\in\N, $ then
\begin{equation}\label{psi'}
\left<\psi'(tu),tu\right>=t^2m_u'(t).
\end{equation}

\noindent {\bf  The case $\Int fu \leq 0$.} 
Notice that the function $m_{u}$
admits an unique turning point $\tilde{t} > 0$, i.e, we get  $m^{\prime}_{u}(t) = 0, t > 0$ if, only if, $t = \tilde{t}$, see Lemma \ref{m_u-comp}.
Moreover, $\tilde{t}$ is a global maximum point for $m_{u}$ such that $m_{u}(\tilde{t}) > 0, m_{u}(\infty) = - \infty$.
As a byproduct there exits an unique $t_1>\tilde t$ such that
$$m_u(t_1)= \Int fu.$$
We emphasize that $m_u'(t_1)<0$, because $m_{u}$ is a decreasing function in $(\tilde{t}, \infty)$. 
Therefore,  using Lemma \ref{m_uegamma_u}, we have  $t_1 u\in \N,$ proving that $\gamma'_u( t_1) =0$. Additionally, by the identity \eqref{psi'}
\begin{equation}\label{rel-m-u-gamma-u}
	m_u(t)=\gamma_u'(t)+ \Int fu,\nonumber
\end{equation}
we get $0>t^2m'_u(t_1)=\left<\psi´(t_1u),t_1u\right>,$
 proving that $t_1u\in\N^-$.

\noindent {\bf  The case $\Int fu > 0$.}  We can consider Lemma \ref{mapf1f2} and we get 
$$m_u(\tilde t) > \Int fu,$$
which  $m_{u}$ is increasing in $(0, \tilde{t})$ and decreasing in $(\tilde{t}, \infty)$.
It is not hard to verify that there exist exactly two points $0<t_1=t_1(u)  <\tilde t < t_2=t_2(u)$ such that
$$m_u(t_i)=\Int fu,~i=1,2,$$
satisfying $m_u'(t_1)>0$ and $m_u'(t_2)<0$. As in the previous step we infer that $t_1u\in\N^+$ and $t_2u\in\N^-$. This completes the proof.  \hfill\cqd

\begin{lem}\label{nehari-} Suppose  either $(f_1)$ or $(f_2)^{\prime}$.
There exist $\delta_1, \lambda_{2} > 0$ in such way that $J(u)\geq\delta_1$ for any $u\in \N^{-}$
where $0<\|f\|_{(\el)'} < \lambda_2$.
\end{lem}
 
 \proof Since $u\in\N^-,$ we have that 
$\left<\psi'(u),u\right><0$. Arguing as in the proof of Lemma \ref{c1}, we obtain 
$$||u||>\left[\Fr{\ell(\ell-1)S^{\frac{\ell^*}{\ell}}}{(\el-1)}\right]^{\frac{1}{\el-\alpha}}.$$

Moreover, in view of \eqref{functional1} and the Sobolev imbedding, we have that
$$
\begin{array}{rcl}
	J(u)&\geq&\ell\left(\Fr{1}{m}-\Fr{1}{\el} \right)\min\{||u||^\ell,||u||^m\}-\left(1-\Fr{1}{\el}\right)\Int fu\\[2ex]
	&\geq&||u||\left[\ell\left(\Fr{1}{m}-\Fr{1}{\el} \right)\min\{||u||^{\ell-1},||u||^{m-1}\}-\left(1-\Fr{1}{\el}\right)S^{-\frac{1}{\ell}}||f||_{(\el)'}\right]. 
\end{array}
$$
By the above inequality, we get 
$$\begin{array}{rcl}J(u)&>&\left[\Fr{\ell(\ell-1)S^{\frac{\ell^*}{\ell}}}{(\el-1)}\right]^{\frac{1}{\el-\alpha}}\left[\ell\left(\Fr{1}{m}-\Fr{1}{\el} \right)\left[\Fr{\ell(\ell-1)S^{\frac{\ell^*}{\ell}}}{(\el-1)}\right]^{\frac{\alpha-1}{\el-\alpha}}-\left(1-\Fr{1}{\el}\right)S^{-\frac{1}{\ell}}||f||_{(\el)'}\right]\\[3ex]
&=&\left[\Fr{\ell(\ell-1)S^{\frac{\ell^*}{\ell}}}{(\el-1)}\right]^{\frac{1}{\el-\alpha}}\left[A-||f||_{(\el)'}B\right].
\end{array}$$ 
Notice that  $A-||f||_{(\el)'}B>0$ if, only if,  $||f||_{(\el)'}<\Fr{A}{B}=\Fr{1}{m}\lambda_1=:\lambda_2,$ where $\lambda_1$ is given by $(f_1).$ On the other hand, if  $(f_2)$ holds, we have
$$||f||_{(\el)'}\leq\min\left\{\frac{1}{m}\lambda_1,\Fr{\el-m}{m-1}\right\}\leq\frac{1}{m}\lambda_1.$$

Hence, in either case  $(f_1)$ ou $(f_2)^{\prime}$, we conclude that  $J(u)\geq \delta_1,$ for all $u\in \mathcal{N}^-$.

\hfill\cqd 

\begin{lem}\label{nehari+}
Suppose $(H)$ and either $(f_1)$ or $(f_2)^{\prime}$. Then, 
$\alpha:=\displaystyle\inf_{u\in \N}J(u) =\alpha^+=\displaystyle\inf_{u\in \N^+} J(u)<0.$ 
\end{lem}

\proof Since $u\in \N^+$ we have that $\left<\psi'(u),u\right>>0$, i.e. $$\Int\phi'(|\nabla u|)|\nabla u|^3+\phi(|\nabla u|)|\nabla u|^2-(\el-1)|u|^{\el}>0.$$
Thus,$$
\begin{array}{rcl}
(\el-1)\Int |u|^{\el}&<&\Int\phi'(|\nabla u|)|\nabla u|^3+\phi(|\nabla u|)|\nabla u|^2
\\[2ex]
&<&(m-1)\Int\phi(|\nabla u|)|\nabla u|^2. \end{array}$$
Consequently,
\begin{equation*}
\int_\Omega |u|^{\el}<\Fr{m-1}{\el-1}\Int\phi(|\nabla u|)|\nabla u|^2
\end{equation*}
On the other hand, if $u\in\N$, using the above inequality and $(\phi_3)$, we get
$$\begin{array}{rcl}
J(u)&\leq& \left(\Fr{1}{\ell}-1\right)\Int\phi(|\nabla u|)|\nabla u|^2+\left(1-\Fr{1}{\el}\right)|u|^{\el}\\[3ex]
&<&\left[\Fr{1-\ell}{\ell }+\Fr{m-1}{\el }\right]\Int\phi(|\nabla u|)|\nabla u|^2<0,
\end{array}$$
  because
	$\left[\Fr{1-\ell}{\ell }+\Fr{m-1}{\el }\right]<0$, since $(H)$ holds. 
Consequently,  \ $\alpha^+<0.$

Since $\N=\N^-\cup\N^+$ and $\alpha^->0$, we have that $\alpha^+=\alpha,$ and the Lemma is proved.
\hfill\cqd

\section{The (PS) condition}

 Here we follow same ideas discussed in Tarantello \cite{tarantello}, in order to prove some auxiliary results  to get the Palais-Smale conditon for the functional $J$ constrained to the Nehari manifold.

\begin{lem}\label{lem1ps} Suppose $(\phi_{1})- (\phi_{3})$ and $(H).$ Let $u \in \mathcal{N}^{+}$ be fixed. Then there exist $\epsilon>0$ and a differentiable function $$\xi:B(0,\epsilon)\subset W_{0}^{1,\Phi}(\Omega)\to (0, \infty),\ \ \ \xi(0)=1,\,\ \xi(v)(u-v)\in \mathcal{N}^{+}, v  \in B(0,\epsilon).$$
Furthermore, we have that
\begin{eqnarray}\label{ps}
\left<\xi'(0),v\right> &=& \dfrac{1}{\left<\psi'(u),u\right>} \Int\left\{[\phi'(|\nabla u|)|\nabla u|+2\phi(|\nabla u|)]\nabla u\nabla v-\el |u|^{\el-2}uv\right.\nonumber\\
&-&\left.fv\right\}.
\end{eqnarray}
\end{lem}
\proof Initially, we define $\psi: \w \backslash \{ 0\} \rightarrow \mathbb{R}$ given by $\psi(u)=\left<J'(u),u\right>$ for $u \in \w \backslash \{ 0\}$. 
Recall that $\left<\psi'(u),u\right> $ is given by \eqref{psi}, and for any $u\in\N,$ $\left<\psi'(u),u\right> $ was defined in Remark \ref{gamma''}.

Now we define $F_u : \mathbb{R} \times \w\backslash \{0\} \rightarrow \mathbb{R}$ given by 
$F_u : \mathbb{R} \times \w\backslash \{0\} \rightarrow \mathbb{R}$ given by $$F_{u}(\xi,w)= \left<J'(\xi(u-w)),\xi(u-w)\right>.$$ Here we observe that $F_u(1,0)=\psi(u)$. As a consequence, for each $u\in\N$, we have
 \begin{eqnarray}
 \partial_1 F_u(1,0)&=&\Int 2\phi(|\nabla u|)|\nabla u|^2+\phi'(|\nabla u|)|\nabla u|^3 \nonumber \\
 &-& \Int \el|u|^{\el}-fu= \left<\psi'(u),u\right>\neq0. \nonumber
 \end{eqnarray}
 
 By using the Inverse Function Theorem, there exist $\epsilon>0$ and a differentiable function $\xi:B(0,\epsilon)\subset W^{1,\Phi}(\Omega)\to (0, \infty)$ satisfying $\xi(0)=1$ and $F_u(\xi(w),w)= \langle J'(\xi(w)(u-w)),\xi(w)(u-w) \rangle = 0,$ i.e. $\xi(w)(u-w)\in\N,\ \ \forall w\in B(0,\epsilon).$ Furthermore, we also have
 \begin{equation*}
 \left<\xi'(w),v\right>=-\Fr{\left<\partial_2F_u(\xi(w),w),v\right>}{\partial_1F_u(\xi(w),w)}, \left<\xi'(0),v\right>=-\Fr{\left<\partial_2F_u(\xi(0),0),v)\right>}{\partial_1F_u(\xi(0),0)}.
 \end{equation*}
 Here $\partial_1F_u$ and $\partial_2F_u$ denote the partial derivatives on the first and second variable, respectively.

 On the other hand, after some manipulations,
  putting $w=0$ and $\xi=\xi(0)=1$, we have 
 $$\begin{array}{rcl}-\left<\partial_2F_u(1,0),v\right>&=&\Int\left(\phi'(|\nabla u|)|\nabla u| +
 2\phi(|\nabla u|)\right)\nabla u\nabla v-\el |u|^{\el-2} uv
 -fv
 \end{array}$$
 Here was used the fact that $\partial_1 F_u(1,0)=\left<\psi'(u),u\right>$ holds for any $u\in\N$. The proof is complete. \hfill\cqd

Similarly, we have the following
\begin{lem}\label{lem2ps}
Suppose $(\phi_{1})- (\phi_{3})$ and $(H)$. Let $u \in \N^{-}$ be fixed. Then there are $\epsilon>0$ and a differentiable function $$\xi^-:B(0,\epsilon)\subset W^{1,\Phi}(\Omega)\to (0, \infty),\ \ \ \xi^-(0)=1,\,\ \xi^-(v)(u-v)\in\N^{-}, \, v \in B(0,\epsilon).$$
Furthermore, we obtain
\begin{eqnarray}\label{ps-}
\left< (\xi^{-})^{\prime}(0),v\right>&=&\dfrac{1}{\left<\psi'(u),u\right>} \Int\left\{[\phi'(|\nabla u|)|\nabla u|+2\phi(|\nabla u|)]\nabla u\nabla v-\el |u|^{\el-2}uv\right.\nonumber\\
&-&\left.fv\right\}.
\end{eqnarray}
\end{lem}

Next, we shall prove that any minimizing sequences on the Nehari manifold in $\mathcal{N}^{+}$ or $\mathcal{N}^{+}$  provides us a Palais-Smale sequence. 

\begin{prop}\label{3.1}Suppose $(\phi_{1})- (\phi_{3})$ and $(H)$.
Then we have the following assertions
\begin{enumerate}
\item there exists a sequence $(u_n)\subset \N$ such that
$J(u_n)= \alpha^++o_n(1)
 \ \mbox{\and } J'(u_n)=o_n(1)\ \mbox{ in}~  W^{-1,\widetilde{\Phi}}(\Omega).$
\item there exists a sequence $(u_n)\subset \N^-$ such that
$J(u_n)=\alpha^-+o_n(1)
 \ \mbox{\and } J'(u_n)=o_n(1)\ \mbox{ in } W^{-1,\widetilde{\Phi}}(\Omega).$
\end{enumerate}
\end{prop}

\begin{prop}
	Suppose $(\phi_{1})- (\phi_{3})$ and $(H)$ hold. Let $(u_n)$ be a minimizing sequence for the functional $J$ constrained to the Nehari mainfold $\mathcal{N}^+$. Then 
	\begin{equation}\label{12b}
	\liminf_{n \rightarrow \infty} ||u_n|| \geq - \alpha^+ \left[\Fr{\el}{(\el-1)||f||_{(\el)'}S^{\frac{-1}{\ell}}}\right] > 0.
	\end{equation}
	and
	\begin{eqnarray}\label{u_n-lim}
	 ||u_n||< \left[\left(\Fr{\el-1}{\el-m}\right)||f||_{(\el)'} S^{\frac{-1}{\ell}} \right]^\frac{1}{\alpha-1},
	\end{eqnarray}
	where $\alpha\in\{\ell,m\}$.  The same property can be proved for the Nehari manifold $\mathcal{N}^-.$ 
\end{prop}
\Proof Remember that $(u_n)\subset\N$, $m\Phi(t)\leq \phi(t)t^2$  and arquing as in the proof of Lemma \ref{nehari-}, we infer that
\begin{equation}
\begin{array}{rcl}\label{12}
0> \alpha^{+} + o_{n}(1) &>&J(u_n)\\ 
&\geq& \Int\left(1-\Fr{m}{\el}\right) \Phi(|\nabla u_n|)-\left(1-\frac{1}{\el}\right)fu
\end{array}
\end{equation}
holds for any $n\in\mathbb{N}$ large enough.  By using the above inequality and the continuous embedding $W_0^{1,\Phi}(\Omega)\hookrightarrow L^{\el}(\Omega)$, we deduce that
\begin{equation*}\label{12a}
||u_n||> \left[-\left(\alpha^++\Fr{1}{n}\right)\Fr{\el}{(\el-1)||f||_{(\el)'}S^{\frac{-1}{\ell}}}\right],
\end{equation*}
and \eqref{12b} holds.

Furthermore, using $\eqref{12}$ and arguing as in (\ref{des2}), we  obtain that
\begin{eqnarray}
\displaystyle \min\{||u_n||^\ell,||u_n||^m\}\leq \Int \Phi(|\nabla u_n|)
&<& \left(\Fr{\el-1}{\el-m}\right)||f||_{(\el)'} S^{\frac{-1}{\ell}} ||u_n||. \nonumber \\
\end{eqnarray}
Hence the last assertions give us
\begin{eqnarray}
||u_n||<\left[\left(\Fr{\el}{\el-m}\right)\left(\Fr{\el -1}{\el }\right)||f||_{(\el)'}S^{\frac{-1}{\ell}}\right]^{\frac{1}{\alpha-1}}= \left[\left(\Fr{\el-1}{\el-m}\right)||f||_{(\el)'} S^{\frac{-1}{\ell}} \right]^\frac{1}{\alpha-1},\nonumber
\end{eqnarray}
where $\alpha\in\{\ell,m\}$.\hfill\cqd

\vskip.2cm

Now we will prove two technical results, which will be used  to prove that any minimizing sequence for $J$ constrained to the Nehari manifold is a Palais-Smale sequence. 

\begin{prop}\label{J'-lim-pont}
	Suppose $(\phi_{1})- (\phi_{3})$ and $(H)$ hold. Then any minimizing sequence $(u_n)$ on the Nehari manifold  $\mathcal{N}^-$ or $\mathcal{N}^+$ satisfies
	\begin{equation}\label{17d}
	\left<J'(u_n), \Fr{u}{||u||}\right>\leq \Fr{C}{n}[||\xi'_n(0)||+1].
	\end{equation}
	where $\xi_n:B_{\frac{1}{n}}(0)\rightarrow\mathbb{R}$ was obtained by Lemmas \ref{lem1ps} and \ref{lem2ps}.
\end{prop}
\Proof Taking $\epsilon_n $ given in Lemma \ref{lem1ps}, put $\rho \in (0, \epsilon_n)$ and $u\in W^{1,\Phi}(\Omega)\backslash \{0\}$. Define the auxiliary function
$$w_\rho=\Fr{\rho u}{||u||}\in B(0,\epsilon_n).$$
Using  Lemma \ref{lem1ps}  we infer that
\begin{equation}\label{177e}
\mu_\rho=\xi(w_\rho)(u_n-w_\rho)\in\N^+ \,\,\mbox{ and } \,\, J(\mu_\rho)-J(u_n)\geq- \frac{1}{n}||\mu_\rho-u_n||.
\end{equation}
Notice also that we have the following convergences 
\begin{equation} \label{17c}
w_\rho\to 0,\ \xi_n(w_\rho)\to 1, \ \mu_\rho\to u_n \mbox{ and } J'(\mu_\rho)\to J'(u_n)
\end{equation}
as $\rho\to 0,$ for any $n \in \mathbb{N}$.

Applying Mean Value Theorem, there exists $t\in(0,1)$ in such way that
$$\begin{array}{rcl}J(\mu_\rho)-J(u_n)
&=&\left<J'(\mu_\rho+t(u_n-\mu_\rho)) -J'(u_n),\mu_\rho-u_n \right>\\[2ex]
&+&\left<J'(u_n),\mu_\rho-u_n\right>.
\end{array}$$
Remind that  $||u_n-\mu_\rho||\to 0$ as $\rho \to 0$. Since  $\mu_\rho\in\N^+$  and using \eqref{177e} and \eqref{17c}, we obtain
\begin{equation*}\label{16}
-\frac{1}{n}||\mu_\rho-u_n|| +o_{\rho}(||\mu_\rho-u_n||)\leq \left<J'(u_n),-w_\rho\right>+(\xi_n(w_\rho)-1)\left<J'(u_n),u_n-w_\rho\right>.
\end{equation*}
where $o_{\rho}(.)$ denotes a quantity that goes to zero as $\rho$ goes to zero.
Using that $\left<J^{\prime}(\mu_\rho),\mu_\rho\right>=0$, we have
\begin{eqnarray*}
	-\frac{1}{n}||\mu_\rho-u_n|| &\leq& o_{\rho}(||\mu_\rho-u_n||) - \rho\left<J'(u_n),\frac{u}{||u||}\right> \nonumber \\
	&+&(\xi_n(w_\rho)-1)\left<J'(u_n)-J'(\mu_\rho),u_n-w_\rho\right>. \nonumber \\
\end{eqnarray*}
From  the above estimates and \eqref{17c} we obtain
\begin{eqnarray}
\left<J'(u_n),\frac{u}{||u||}\right> &\leq& \Fr{||\mu_\rho-u_n||}{n\rho}+\Fr{o_{\rho}(||\mu_\rho-u_n||)}{\rho}
\nonumber \\
&+& \Fr{(\xi_n(w_\rho)-1)}{\rho}\left<J'(u_n)-J'(\mu_\rho),u_n-w_\rho\right>.  \nonumber
\end{eqnarray}
Noticing that

$$\ds\lim_{\rho\to 0}\Fr{|\xi_n(w_\rho)-1|}{\rho}=\left<\xi'_n(0),\frac{u}{||u||}\right>\leq ||\xi'_n(0)||,$$
from this inequality we have

\begin{equation}\label{17a}
||\mu_\rho-u_n||\leq \rho|\xi_n(w_\rho)|+|\xi_n(w_\rho)-1| \,||u_n||\mbox { and } \ds\lim_{\rho\to 0}\Fr{|\xi_n(w_\rho)-1|}{\rho}\leq||\xi'_n(0)||.
\end{equation}

Therefore, using the fact that $(u_n)$ is bounded and \eqref{17a}, we infer that
$$\begin{array}{rcl}
\ds\lim_{\rho\to 0} \Fr{||\mu_\rho-u_n||}{\rho n}&\leq&\ds\lim_{\rho\to 0}\Fr{1}{n}\left[||\xi_n(w_\rho)||+\Fr{|\xi_n(w_\rho)-1|}{\rho}||u_n||\right]\\[2ex]
&\leq& \Fr{1}{n}\left[1+||\xi'_n(0)||\,\,||u_n||\right]\leq \Fr{C}{n}\left[1+||\xi'_n(0)||\right].
\end{array}$$

On the other hand, since $\Fr{\xi_n(w_\rho)-1}{\rho}$ and $\xi_n(w_\rho)$ are bounded for $\rho>0$ small enough, we obtain
$$\begin{array}{rcl}
\|\mu_\rho-u_n \|&=&|\rho|\left|\left|\Fr{\xi_n(w_\rho)-1}{\rho}u_n-\xi_n(w_\rho)\Fr{u}{||u||}\right|\right|\\[2ex]
&\leq& |\rho|\left[\left|\Fr{\xi_n(w_\rho)-1}{\rho}\right|||u_n||+|\xi_n(w_\rho)|\right].
\end{array}$$
Since $(u_n)$ is bounded there exists a constant $C>0$ in such  that
\begin{equation*}\label{17b}
\Fr{||\mu_\rho-u_n||}{\rho}\leq C[||\xi'_n(0)||+1].
\end{equation*}
Putting all these estimates together we prove \eqref{17d} holds. 
\hfill\cqd

\begin{prop}\label{xi_n-bound}
	Under the hypotheses of Proposition \ref{J'-lim-pont} there exists $C>0$ such that 
	$$||\xi'_n(0)||\leq C,~\forall n\in\mathbb{N}.$$
\end{prop}
\Proof Firstly notice that the numerator in \eqref{ps} is bounded from below away from zero by $ b ||v||$ where $b> 0$ is a constant. Define the auxiliary function $\chi_n : \w \rightarrow \mathbb{R}$ given by
$$\begin{array}{rcl}\chi_n(v)&=&\Int[\phi'(|\nabla u_n|)|\nabla u_n|+2\phi(|\nabla u_n|)]\nabla u_n\nabla v\\[2ex]&-&\Int\el |u_n|^{\el-2}u_nv-fv.\end{array}$$
Using  that $\Fr{|\phi'(t)t|}{\phi(t)}\leq \max\{|\ell-2|,|m-2|\}:=C_1$ and Holder's inequality, we obtain
$$\begin{array}{rcl}|\chi_n(v)|&\leq& C_1\Int\phi(|\nabla u_n|)|\nabla u_n||\nabla v|
+\el\Int |u_n|^{\el-1}|v|+\Int |f||v|\\[2ex]
&\leq& 2C_1||\phi(|\nabla u_n|)|\nabla u_n|||_{\tilde{\Phi}}||v||
+\el\Int |u_n|^{\el-1}|v|+\Int|f||v|\\[2ex]
&\leq& C_2 \max\left\{\left(\Int \tilde{\Phi}(\phi(|\nabla u_n|))|\nabla u_n|\right)^{\frac{\ell-1}{\ell}},\left(\Int \tilde{\Phi}(\phi(|\nabla u_n|))|\nabla u_n|\right)^{\frac{m-1}{m}}\right\} ||v|| \\[2ex]
&+&\el\Int |u_n|^{\el-1}|v|+\Int |f||v|.
\end{array}
$$
In virtue of  the  inequality $\widetilde\Phi(t\phi(t))\leq \Phi(2t)\leq 2^m\Phi(t), t \geq 0$ and \eqref{u_n-lim} there exists a constant $C_3>0$ such that
$$\begin{array}{rcl}|\chi_n(v)|
&\leq& C_3\max\left\{\left(\Int \Phi(|\nabla u_n|)\right)^{\frac{\ell-1}{\ell}},\left(\Int \Phi(|\nabla u_n|)\right)^{\frac{m-1}{m}}\right\}||v||\\[2ex]
&+&\el\Int |u_n|^{\el-1}|v|+\Int |f||v|\\[2ex]
&\leq&
C_3 || u_n||^\beta||v||+\el\Int |u_n|^{\el-1}|v|+\Int |f||v|\\[2ex]
&\leq& C_4||v||+\el\Int |u_n|^{\el-1}|v|+\Int |f||v|.
\end{array}
$$ where $\beta\in\{\ell-1,\frac{\ell}{m}(\ell-1),m-1,\frac{m}{\ell}(m-1)\}$.

We shall estimate the terms $\Int |u_n|^{\el-1}|v| $ and $\Int |f||v|$.
Employing Holder's inequality and Sobolev imbedding we obtain
$$\Int |u_n|^{\el-1}|v|\leq \left(\Int|u_n|^{\el}\right)^{\frac{\el-1}{\el}}\left(\Int|v|^{\el}\right)^{\frac{1}{\el}}\leq C_5 \|u_n\|^{\ell^{*}-1} \|v\| \leq C_6 \|v\|.$$ and 
$$\Int |f||v|\leq \left(\Int|f|^{(\el)'}\right)^{\frac{1}{(\el)'}}\left(\Int|v|^{\el}\right)^{\frac{1}{\el}}=||f||_{(\el)'}||v||_{\el}\leq S^{\frac{-1}{\ell}}||f||_{(\el)'}||v||.$$ 

Combining  the estimates above there exists a constant $c>0$ in such that $|\chi_n(v)|\leq c||v||$. 

Next, we will  show that there exists a constant $d>0$, independent in $n$, such  that $\gamma''_{u_n}(1)\geq d$. Indeed, arguing by contradiction that $\gamma''_{u_n}(1) = o_n(1)$. It follows from \eqref{12b} that there exists $a_\lambda>0$ satisfying
\begin{equation}\label{nozero}
\ds\liminf_{n \rightarrow \infty}||u_n||\geq a>0
\end{equation}

 Using \eqref{nehari} and \eqref{eq5}, as well as,  $\left<\phi'(u_n),u_n\right> = o_n(1),$ we deduce that
$$o_n(1)=\left<\phi'(u_n),u_n\right>=\Int \phi(|\nabla u_n| )|\nabla u_n|^2+\phi'(|\nabla u_n|)||\nabla u|^3+(1-\el)|u_n|^{\el}.$$
Under hypothesis $(\phi_3)$ and the  Sobolev embeddings we infer that
\begin{eqnarray}
(\ell-1)\Int \phi(|\nabla u| )|\nabla u_n|^2
&\leq& (\el-1)S^{\frac{-\el}{\ell}}||u_n||^{\el}+o_n(1).\nonumber
\end{eqnarray}

On the other hand, we observe that
$$(\ell-1)\Int \phi(|\nabla u_n| )|\nabla u_n|^2dx\geq\ell(\ell-1)\Int \Phi(|\nabla u_n|)\geq \ell(\ell-1)\min\{||u_n||^\ell,||u_n||^m\}. $$
Using the above estimates we get
$$\ell(\ell-1)\min\{||u_n||^\ell,||u_n||^m\}\leq(\el-1)S^{\frac{-\el}{\ell}}||u_n||^{\el}+o_n(1).$$
Hence, we have 
$$\ell(\ell-1)\leq {(\el-1)S^{\frac{-\el}{\ell}}}||u_n||^{\el-\alpha}+\displaystyle\frac{o_n(1)}{||u_n||^\alpha}$$
where $\alpha=\ell$ whenever $||u_n||\geq 1$ and $\alpha=m$ whenever $||u_n||\leq 1$. Furthermore, using \eqref{nozero}, we obtain
\begin{equation}\label{des1aa}
||u_n||\geq \left[\Fr{\ell(\ell-1)}{(\el-q)S^{\frac{-\el}{\ell}}}\right]^{\frac{1}{\el-\alpha}}+o_n(1).
\end{equation}
Using  \eqref{eq5}, $(\phi_3)$ and Holder inequality, we obtain
\begin{eqnarray}
(\el-m)\Int \phi(|\nabla u| )|\nabla u_n|^2
&\leq&(\el-1)S^{\frac{-1}{\ell}}||f||_{(\el)'}||u_n||+o_n(1).\nonumber
\end{eqnarray}
Combining the above inequalities, we get
$$\Fr{\ell(\el-m)}{(\el-1)S^{\frac{-1}{\el}}||f||_{(\el)'}}||u_n||^\alpha=\Fr{\ell(\el-m)}{(\el-1)S^{\frac{-1}{\ell}}||f||_{({\el)'}}}\min\{||u_n||^\ell,||u_n||^m\}\leq ||u_n||+o_n(1).$$
To sum up, using the estimate \eqref{nozero}, we can be shown that
\begin{eqnarray}
||u_n||
&\leq&\left[ \Fr{(\el-1)S^{\frac{-1}{\ell}}\|f\|_{\left(\el\right)'}}{\ell(\el-m)}\right]^{\frac{1}{\alpha-1}}+o_n(1).\nonumber
\end{eqnarray}
Arguing as in the proof of Lemma \ref{c1}, by the above inequality and \eqref{des1aa}  we have a  contradiction since either $(f_1)$ or $(f_2)$ hold.  This completes the proof. \hfill\cqd
\vskip.2cm
\nd{\bf Proof of Proposition \ref{3.1}}
We shall prove the item $(1)$. The proof of item $(2)$ follows similarly using Lemma \ref{lem2ps} instead of Lemma \ref{lem1ps}.  Applying Ekeland's variational principle there exists a sequence $(u_n)\subset\N^+$ in such  way that
\begin{description}
\item[(i)] $J(u_n) \,= \,\alpha^+ + o_n(1)$,
\item [(ii)] $J(u_n)<J(w)+\frac{1}{n}||w-u||, \,\,\forall \,\, w\,\in\N^+.$
\end{description}

In what follows we shall prove that $\displaystyle \lim_{n \rightarrow \infty} ||J'(u_n)||\to 0$. From Proposition \ref{xi_n-bound}, there exist $C>0$  independent on $n\in\mathbb{N}$ such that $\|\xi_n(0)\|\leq C$.  This estimate together with Proposition \ref{J'-lim-pont}
$$
\left<J'(u_n), \Fr{u}{||u||}\right>\leq \Fr{C}{n},~u\in\w/\{0\}.
$$
This implies that  $\|J'(u_n)\|\rightarrow 0$ as $n\rightarrow\infty$. This finishes the proof. \hfill\cqd


\section{The proof of our main theorems}

\subsection{The proof of Theorem \ref{teorem1}}
We are going to apply the following result, whose proof is made by using the concentration compactness principle due to Lions for Orlicz -Sobolev framework, see \cite{Willem} or else in \cite{CSGG,Fuk_1}.
\begin{lem}\label{conv_grad_qtp}
	$(i)$ $\phi(|\nabla u_n|)\nabla u_n\rightharpoonup \phi(|\nabla u|)\nabla u$ in $\prod L_{\widetilde{\Phi}}(\Omega)$;\\
	$(ii)$ $|u_n|^{\ell^*-2}u_n \rightharpoonup|u|^{\ell^*-2}u$ in $L^{\frac{\ell^*}{\ell^*-1}}(\Omega)$.
\end{lem}
Let $\|f\|_{(\el)'} < \Lambda_1 = \min\left\{\lambda_1,\Fr{\el-m}{m-1}\right\}$  where $\lambda_1 > 0$ is given by $(f_1)$. 

From Lemma \ref{nehari+} we infer that $$\alpha^+:=\ds\inf_{u\in\N^+}J(u)=\ds\inf_{u\in\N}J(u) < 0.$$
We will find a function $u\in \N^+$ in such that $$J(u)= \ds\min_{u\in \N^+}J (u)=:\alpha^+ \,\, \mbox{and} \,\, J^{\prime}(u) \equiv 0.$$
First of all, using Proposition \ref{lem1ps}, there exists a minimizing sequence denoted by $(u_n)\subset W^{1,\Phi}(\Omega)$ such that
\begin{equation}\label{cerami1}
J(u_n)=\alpha^++o_n(1) \mbox{ and }
J'(u_n)=o_n(1).
\end{equation}
Since the functional $J$ is coercive in $\N^+$, this implies   that $(u_n)$  is bounded in $\N^+$. Therefore, there exists a function $u\in{W^{1,\Phi}_0(\Omega)}$   such  that
\begin{equation} \label{convergencia}
u_n \rightharpoonup u \,\, \mbox{ in } \,\, W_0^{1,\Phi}(\Omega),~~
u_n \to u \,\,\mbox{a.e.}\,\, \mbox{ in } \Omega,~~
u_n \to u \,\, \mbox{ in } \,\, L^{\Phi}(\Omega).
\end{equation}
\nd We shall prove that $u$ is a weak solution for the problem elliptic problem \eqref{eq1}.
Notice that, by \eqref{cerami1}, we mention that
$$\on=\left<J'(u_n),v\right>=\Int \phi(|\nabla u_n|)\nabla u_n\nabla v- fv -|u_n|^{\el-2}u_nv$$
holds for any $v \in \w$. In view of \eqref{convergencia} and Lemma \ref{conv_grad_qtp} we get
$$\Int \phi(|\nabla u|)\nabla u\nabla v-f v-|u|^{\el-2}uv = 0$$
for any $v\in W^{1,\Phi}(\Omega)$ proving that u is a weak solution to the elliptic problem \eqref{eq1}. In addition, the weak solution $u$ is not zero. In fact, using the fact that $u_n\in \mathcal{N}^{+},$ we obtain
$$\begin{array}{rcl}
\Int fu_n&=& \Int( \Phi(|\nabla u_n|)- \Fr{1}{\el}\phi(|\nabla u_n|)|\nabla u_n|^2) \Fr{\el}{\el-1} -J(u_n) \Fr{\el}{\el-1}\\[3ex]
&\geq& \Fr{\el}{\el-1} \left(1 - \dfrac{m}{\ell^*}\right)\Int \Phi(|\nabla u_n|) -J(u_n)\Fr{\el}{\el-1} \\[3ex]
&\geq& -J(u_n)\Fr{\el}{\el-1}.
\end{array}
$$
From \eqref{cerami1} and \eqref{convergencia} we obtain 
\begin{eqnarray}\label{fu-pos}
	\Int fu\geq-\alpha^{+}\Fr{\el}{\el-1} > 0.
\end{eqnarray}
Hence $u\not\equiv 0$.

We shall prove that $J(u)=\alpha^+$ and $u_n\to u$ in $W_0^{1,\Phi}(\Omega)$.
Since $u\in \N$ we also see that
$$\alpha^+ \leq J(u)=\Int \Phi(|\nabla u|)-\Fr{1}{\el}\phi(|\nabla u|)|\nabla u|^2-\left(1-\Fr{1}{\el}\right)fu.$$
Notice that
$$t\mapsto\Phi(t)-\Fr{1}{\el}\phi(t)t^2$$
is a convex function. In fact, by hypothesis $(\phi_3)$ and $m<\el$, we infer that
\begin{eqnarray}
	\left(\Phi(t)-\Fr{1}{\el}\phi(t)t^2\right)''&=&\left[ \left(1-\frac{1}{\el}\right)t\phi(t)-\frac{1}{\el}t(t\phi(t))'\right]'\nonumber\\
	&=& (t\phi(t))' \left[\left(1-\frac{2}{\el}\right)-\frac{1}{\el}\frac{t(t\phi(t))''}{(t\phi(t))'}\right]\nonumber\\
	&\geq& (t\phi(t))'\left(1-\frac{m}{\el}\right)>0, t > 0.\nonumber
\end{eqnarray}
In addition, the last assertion says that
$$u\longmapsto \Int \Phi(|\nabla u |)-\Fr{1}{\el}\phi(|\nabla u |)|\nabla u |^2dx$$
is weakly lower semicontinuous function. Therefore we obtain
\begin{eqnarray}
	\alpha^+ \leq J(u) &\leq & \liminf \left(\Int \Phi(|\nabla u_n|)-\Fr{1}{\el}\phi(|\nabla u_n|)|\nabla u_n|^2
	-\left(1-\Fr{1}{\el}\right)fu_n\right)\nonumber\\
	&=&\liminf J(u_n)= \alpha^+.\nonumber
\end{eqnarray}
This implies that  $J(u)=\alpha^+.$ Additionally, using \eqref{convergencia}, we also have
$$\begin{array}{rcl}
J(u)&=&\Int \Phi(|\nabla u|)-\Fr{1}{\el}\phi(|\nabla u|)|\nabla u|^2-\left(1-\Fr{1}{\el}\right)fu\\[3ex]
&=&
\lim \left(\Int \Phi(|\nabla u_n|)-\Fr{1}{\el}\phi(|\nabla u_n|)|\nabla u_n|^2\right)-\left(1-\Fr{1}{\el}\right)\Int fu.
\end{array}$$
From the last identity 
$$\lim \left(\Int \Phi(|\nabla u_n|)-\Fr{1}{\el}\phi(|\nabla u_n|)|\nabla u_n|^2\right)=\Int \Phi(|\nabla u|)-\Fr{1}{\el}\phi(|\nabla u|)|\nabla u|^2.$$
In view of Brezis-Lieb Lemma, choosing $v_n=u_n-u,$ we infer that
\begin{eqnarray}
	\lim \left(\Int \Phi(|\nabla u_n|)-\Fr{1}{\el}\phi(|\nabla u_n|)|\nabla u_n|^2+ \Phi(|\nabla v_n|)-\Fr{1}{\el}\phi(|\nabla v_n|)|\nabla v_n|^2\right)\nonumber\\
	=\Int \Phi(|\nabla u|)-\Fr{1}{\el}\phi(|\nabla u|)|\nabla u|^2.
\end{eqnarray}
 The previous assertion implies that
$$0=\lim \left(\Int \Phi(|\nabla v_n|)-\Fr{1}{\el}\phi(|\nabla v_n|)|\nabla v_n|^2\right)\geq \lim\left(1-\Fr{m}{\el}\right)\Int \Phi(|\nabla v_n|)\geq 0.$$ Therefore, we obtain that $\lim \int_{\Omega} \Phi(|\nabla v_n|)=0$ and $u_n\to u \,\,\mbox{in} \,\, W^{1,\Phi}(\Omega).$  Hence we conclude that $u_{n} \rightarrow u$ in $\w$.

We shall prove that $u\in \N^+$. Arguing by contradiction we have that $u\notin\N^+$. Using Lemma \ref{fib} there are unique $t_0^+, t_0^->0$ in such way that $t_0^+u\in\N^+$ and $t_0^-u\in\N^-$. In particular, we know that $t_0^+< t_0^-=1.$ Since
$$\Fr{d}{dt}J(t_0^+u)=0$$ 
and using \eqref{fu-pos} together the Lemma \ref{fib} we have that
$$\Fr{d}{dt}J(tu)>0,~t\in (t_0^+, t_0^-).$$
So, there exist $t^- \in ( t_0^+, t_0^- )$ such that $J(t_0^+u)<J(t^-u)$.

\nd In addition  $J(t_0^+u)<J(t^-u)\leq J(t_0^-u)=J(u)$ which is a contradiction to the fact that $u$ is a minimizer in $\N^+$. So that $u$ is in $\N^+$.

\nd To conclude the proof of  theorem it remains to show that $ u\geq 0 $ when $ f \geq0.$ For this we will argue as in \cite{tarantello}. Since $u \in \N^+$, by Lemma \ref{fib} there exists a $ t_0 \geq 1 $ such that $ t_0 | u | \in \N^+ $ and $t_0|u|\geq |u|.$ Therefore if $f\geq 0$, we get
	$$J(u)=\displaystyle\inf_{w\in\N^+}J(w)\leq J(t_0|u|)\leq J(|u|)\leq J(u).$$
\nd So we can assume without loss of generality that $u\geq 0.$

\subsection{The proof of Theorem \ref{teorema2}}
Let $||f||_{(\el)'} < \Lambda_2 = \min\left\{\lambda_2,\Fr{\el-m}{m-1}\right\}$  where $\lambda_2 > 0$ is given by Lemma \ref{nehari-}.

First of all, from Lemma \ref{nehari-}, there exists $\delta_1>0$ such that $J(v)\geq \delta_1$ for any $v\in \N^{-}.$
So that,
$$\alpha^{-}:= \ds \inf_{v \in \N^{-}}J(v)\geq \delta_1>0.$$

Now we shall consider a minimizing sequence $(v_n)\subset \N^{-}$ given in Proposition \ref{lem1ps}, i.e, $(v_n)\subset \N^{-}$ is a sequence satisfying
\begin{equation} \label{e1}
\ds\lim_{n\to\infty}J(v_n)=\alpha^{-} \,\,\mbox{and} \,\, \ds\lim_{n\to\infty} J^{\prime}(v_{n}) = 0.
\end{equation}
Since $J$ is coercive in $\N$ and so on $\N^{-}$, using Lemma \ref{c1}, we have that $(v_n)$ is  bounded sequence
in $W^{1,\Phi}_{0}(\Omega).$ Up to a subsequence we assume that $v_n\rightharpoonup v$ in $W^{1,\Phi}_{0}(\Omega)$ holds for some $v \in \w$.  Additionally, using the fact that  $\ell^{*}>1$, we get $t<<\Phi_{*}(t)$ and $W_{0}^{1,\Phi}(\Omega)\hookrightarrow L^1(\Omega)$ is also a compact embedding. This fact implies that $v_n\to v$ in  $L^{1}(\Omega).$ In this way, we can obtain
\begin{equation*}\label{lim1}
\ds\lim_{n\to\infty} \Int fv_n=\Int fv.
\end{equation*}

Now we claim that $v \in \w$ given just above is a weak solution to the elliptic problem \eqref{eq1}. In fact, using \eqref{e1}, we infer that
$$\left<J'(v_n), w\right>=\Int \phi(|\nabla v_n|)\nabla v_n\nabla w-fw-|v_n|^{\el-2}v_n w = o_{n}(1)$$
holds for any $w \in \w$. Now using Lemma \ref{conv_grad_qtp} we get
$$\Int \phi(|\nabla v|)\nabla v\nabla w-fw -|v|^{\el-2}v w = 0, w \in \w.$$
So that $v$ is a critical point for the functional $J$.
Without any loss of generality, changing the sequence $(v_{n})$ by $(|v_{n}|)$, we can assume that $v \geq 0$ in $\Omega$.

Next we claim that $v \neq 0$. The proof for this claim follows arguing by contradiction assuming that $v \equiv 0$. Recall that
$J(t v_{n}) \leq J(v_{n})$ for any $t \geq 0$ and $n \in \mathbb{N}$. These facts together with Lemma \ref{lema_naru} imply that
\begin{eqnarray}
\left(1 - \dfrac{m}{\ell^{*}}\right)\int_{\Omega} \Phi(|\nabla t v_{n}|) &\leq&  \left(t - 1\right)\left(1-\Fr{1}{\el}\right) \int_{\Omega} f v_n\nonumber \\
&+& \left(1 - \dfrac{\ell}{\ell^{*}}\right)\int_{\Omega} \Phi(|\nabla v_{n}|). \nonumber
\end{eqnarray}
Using the above estimate, Lemma \ref{lema_naru} and the fact that $(v_{n})$ is bounded, we obtain
\begin{equation*}
\min(t^{\ell}, t^{m}) \left(1 - \dfrac{m}{\ell^{*}}\right) \int_{\Omega} \Phi(|\nabla v_{n}|) \leq \left(t - 1\right)\left(1-\Fr{1}{\el}\right) \int_{\Omega} fv_n + C
\end{equation*}
holds for some $C > 0$. These inequalities give us
\begin{equation*}
\begin{array}{rcl}\min(t^{\ell}, t^{m}) \left(1 - \dfrac{m}{\ell^{*}}\right) \Int \Phi(|\nabla v_{n}|) 
&\leq&  \left(t - 1\right)\left(1-\Fr{1}{\el}\right)S^{\frac{-1}{\ell}} ||f||_{(\el)'}||v_n|| + C.\end{array}
\end{equation*}

It is no hard to verify that  $\|v_{n}\| \geq c > 0$ for any $n \in \mathbb{N}$. Using  Proposition \ref{lema_naru} we get
\begin{equation*}
 \min(t^{\ell}, t^{m}) \leq o_{n}(1) t + C
\end{equation*}
holds for any $t \geq 0$ where $C = C(\ell,m,\ell^{*}, \Omega, a,b) > 0$ where $o_{n}(1)$ denotes a quantity that goes to zero as $n \rightarrow \infty$. Here was used the fact $v_{n} \rightarrow 0$ in $L^1(\Omega)$. This estimate does not make sense for any $t > 0$ big enough. Hence $v \neq 0$ as claimed. Hence $v$ is in $\mathcal{N} = \mathcal{N^{+}} \cup \mathcal{N^{+}}$.

Next, we shall prove that $v_n\to v$ in $W_0^{1,\Phi}(\Omega)$. The proof follows arguing by contradiction.
Assume that $\displaystyle \liminf_{n \rightarrow \infty} \int_{\Omega} \Phi(|\nabla v_{n} - \nabla v|) \geq \delta$ holds for some $\delta > 0$.
Recall that $\Psi: \mathbb{R} \rightarrow \mathbb{R}$ given by
$$t\mapsto \Psi(t) := \Phi(t)-\Fr{1}{\el}\phi(t)t^2$$
is a convex function for each $t \geq 0$. The Brezis-Lieb Lemma for convex functions says that
\begin{equation*}
\lim_{n \rightarrow \infty} \int_{\Omega} \Psi(|\nabla v_{n}|) - \Psi(|\nabla v_{n}- v|) =  \int_{\Omega} \Psi(|\nabla v|)
\end{equation*}
In particular, the last estimate give us
\begin{equation*}
\int_{\Omega} \Psi(|\nabla v|) < \liminf_{n \rightarrow \infty}  \int_{\Omega} \Psi(|\nabla v_{n}|).
\end{equation*}
Since $v\in \N$ there exists unique $t_{0}$ in $(0, \infty)$ such that $t_{0} v \in \mathcal{N}^{-}$. It is easy to verify that
\begin{equation*}
\int_{\Omega} \Psi(|\nabla t_{0}v|) < \liminf_{n \rightarrow \infty}  \int_{\Omega} \Psi(|\nabla t_{0} v_{n}|).
\end{equation*}
This implies that
\begin{eqnarray}
\alpha^{-}&\leq& J(t_{0}v ) = \Int \Psi(|\nabla t_{0} v|)-\left(1-\Fr{1}{\el}\right)t_{0}f \nonumber \\
 &<& \liminf_{n \rightarrow \infty} \Int \Psi(|\nabla t_{0} v_{n}|)-\left({1}-\Fr{1}{\el}\right)t_{0} fv_n \nonumber \\
 &=& \liminf_{n \rightarrow \infty} J(t_{0}v_{n}) \leq \liminf_{n \rightarrow \infty} J(v_{n}) = \alpha^{-}.  \nonumber
\end{eqnarray}
This is a contradiction proving that $v_n\to v$ in $W_0^{1,\Phi}(\Omega)$. Therefore $v$ is in $\mathcal{N}^{-}$. This follows from the strong convergence and the fact that $t = 1$ is the unique maximum point for the fibering map $\gamma_{v}$ for any $v \in \mathcal{N}^{-}$.
Hence using the same ideas discussed in the proof of Theorem \ref{teorem1} we infer that
\begin{equation*}
\alpha^{-} \leq J(v) \leq \liminf J(v_{n}) = \alpha^{-}.
\end{equation*}
In particular, $\alpha^{-} = J(v)$ and
\begin{equation*}
\lim \int_{\Omega} \Phi(|\nabla v_{n}|) - \dfrac{1}{\ell^{*}} \phi(|\nabla v_{n}|)|\nabla v_{n}|^{2} = \int_{\Omega} \Phi(|\nabla v|) - \dfrac{1}{\ell^{*}} \phi(|\nabla v|)|\nabla v|^{2}.
\end{equation*}
Hence,  $J(v)\geq \delta_{1} > 0$. This finishes the proof of Theorem \ref{teorema2}.

 \hfill\cqd

\subsection{The proof of Theorem \ref{teorema3}}

In view of Theorems \ref{teorem1} and \ref{teorema2} there are $u\in \mathcal{N}^{+}$ and $v \in \mathcal{N}^{-}$ in such way that $$J(u)=\ds\inf_{w\in \mathcal{N}^{+}}J(w)\ \ \ \mbox{and}\ \ \  J(v)=\ds\inf_{ w\in \mathcal{N}^{-}}J(w).$$
Using  that $0 < ||f||_{(\el)'} < \Lambda := \min\{\Lambda_1,\Lambda_2\}$ where $\Lambda_1, \Lambda_2>0$ are given by Theorem \ref{teorem1} and Theorem \ref{teorema2} we stress that $\N^{+}\cap \N^{-}= \emptyset$.

Therefore, $u,v$ are nonnegative  solutions to the elliptic problem \eqref{eq1}, ($u$ being a ground state solution), whenever $0 < ||f||_{(\el)} < \Lambda$.
This completes the proof.
\hfill\cqd

\section{Appendix}

The reader is  referred to  \cite{A,Rao1} regarding Orlicz-Sobolev spaces.  The usual norm on $L_{\Phi}(\Omega)$ is ( Luxemburg norm),
$$
\|u\|_\Phi=\inf\left\{\lambda>0~|~\int_\Omega \Phi\left(\frac{u(x)}{\lambda}\right) dx \leq 1\right\}
$$
and  the  Orlicz-Sobolev norm of $ W^{1, \Phi}(\Omega)$ is
\[
  \displaystyle \|u\|_{1,\Phi}=\|u\|_\Phi+\sum_{i=1}^N\left\|\frac{\partial u}{\partial x_i}\right\|_\Phi.
\]
Recall that
$$
\widetilde{\Phi}(t) = \displaystyle \max_{s \geq 0} \{ts - \Phi(s) \},~ t \geq 0.
$$
\nd It turns out that  $\Phi$ and $\widetilde{\Phi}$  are  N-functions  satisfying  the $\Delta_2$-condition, (cf. \cite[p 22]{Rao1}).
In addition,   $L_{\Phi}(\Omega)$  and $W^{1,\Phi}(\Omega)$  are separable, reflexive,  Banach spaces.

Using the Poincar\'e inequality for the $\Phi$-Laplacian operator  it follows that
\[
\|u\|_\Phi\leq C \|\nabla u\|_\Phi~\mbox{for any}~ u \in W_{0}^{1,\Phi}(\Omega)
\]
holds true for some $C > 0$, see Gossez \cite{Gz1,gossez-Czech}.
\nd As a consequence,  $\|u\| :=\|\nabla u\|_\Phi$ defines a norm in $W_{0}^{1,\Phi}(\Omega)$, equivalent to $\|.\|_{1,\Phi}$. Let $\Phi_*$ be the inverse of the function
$$
t\in(0,\infty)\mapsto\int_0^t\frac{\Phi^{-1}(s)}{s^{\frac{N+1}{N}}}ds
$$
\nd which extends to ${\mathbb{R}}$ by  $\Phi_*(t)=\Phi_*(-t)$ for  $t\leq 0.$
We say that a N-function $\Psi$ grow essentially more slowly than $\Phi_*$, we write $\Psi<<\Phi_*$, if
$$
\lim_{t\rightarrow \infty}\frac{\Psi(\lambda t)}{\Phi_*(t)}=0,~~\mbox{for all}~~\lambda >0.
$$

The compact embedding below (cf. \cite{A, DT}) will be used in this paper:
$$
\displaystyle W_{0}^{1,\Phi}(\Omega) \stackrel{\tiny cpt}\hookrightarrow L_\Psi(\Omega),~~\mbox{if}~~\Psi<<\Phi_*,
$$
in particular, as $\Phi<<\Phi_*$ (cf. \cite[Lemma 4.14]{Gz1}),
$$
W_{0}^{1,\Phi}(\Omega) \stackrel{\tiny{cpt}} \hookrightarrow L_\Phi(\Omega),
$$
 
Furthermore, the following continuous embeddings hold (see \cite{A, DT,Gz1})
$$
W_0^{1,\Phi}(\Omega) \stackrel{\mbox{\tiny cont}}{\hookrightarrow} L_{\Phi_*}(\Omega),
$$
 $$ L_\Phi(\Omega)\overset{\mbox{\tiny cont}}{\hookrightarrow}  L^{\ell}(\Omega)\ \ \mbox{and}\ \ L_{\Phi_{*}}(\Omega)\overset{\mbox{\tiny cont}}{\hookrightarrow} L^{\ell^{*}}(\Omega).$$

\begin{rmk}\label{rmk-psi}
    The function $\psi(t) = t^{r-1}, r \in [1, \ell^*)$ satisfies $\Psi<<\Phi_*$ where $\Psi(t) = \int_{0}^{t} \psi(s) ds, t \in \mathbb{R}$. In other words,  the function $\Psi$ grow essentially more slowly than $\Phi_*$. In fact, we easily see that
    $$\lim_{t\rightarrow\infty}\frac{\Psi(\lambda t)}{\Phi_*(t)}\leq \frac{\lambda^{r}}{r\Phi_*(1)}\lim_{t\rightarrow\infty}\frac{1}{t^{\ell^*-r}}=0,~~\mbox{for all}~~\lambda>0.$$
    In that case $W_{0}^{1,\Phi}(\Omega) \stackrel{cpt}\hookrightarrow L_\Psi(\Omega)$.
\end{rmk}

Now we refer the reader to \cite{Fuk_1,fang} for some  elementary results on Orlicz and Orlicz-Sobolev spaces.
\vskip.2cm

\begin{prop}\label{lema_naru}
       Assume that  $\phi$ satisfies  $(\phi_1)-(\phi_3)$.
        Set
         $$
         \zeta_0(t)=\min\{t^\ell,t^m\},~~~ \zeta_1(t)=\max\{t^\ell,t^m\},~~ t\geq 0.
        $$
  \nd Then  $\Phi$ satisfies
       $$
            \zeta_0(t)\Phi(\rho)\leq\Phi(\rho t)\leq \zeta_1(t)\Phi(\rho),~~ \rho, t> 0,
        $$
$$
\zeta_0(\|u\|_{\Phi})\leq\int_\Omega\Phi(u)dx\leq \zeta_1(\|u\|_{\Phi}),~ u\in L_{\Phi}(\Omega).
 $$
    \end{prop}

\begin{prop}\label{fang}
       Assume that $(\phi_1)-(\phi_3)$ holds.
        Define the function
         $$
         \eta_0(t)=\min\{t^{\ell-2},t^{m -2}\},~~~ \eta_1(t)=\max\{t^{\ell -2},t^{m - 2}\},~~ t\geq 0.
        $$
  \nd Then the function $\phi$ verifies
       $$
            \eta_0(t)\phi(\rho) \leq \phi(\rho t)\leq \eta_1(t)\phi(\rho),~~ \rho, t> 0,
        $$
    \end{prop}

\begin{prop}\label{lema_naru_*}
    Assume that  $\phi$ satisfies $(\phi_1)-(\phi_3)$.  Set
    $$
    \zeta_2(t)=\min\{t^{\ell^*},t^{m^*}\},~~ \zeta_3(t)=\max\{t^{\ell^*},t^{m^*}\},~~  t\geq 0
    $$
\nd where $1<\ell,m<N$ and $m^* = \frac{mN}{N-m}$, $\ell^* = \frac{\ell N}{N-\ell}$.  Then
        $$
            \ell^*\leq\frac{t^2\Phi'_*(t)}{\Phi_*(t)}\leq m^*,~t>0,
       $$
        $$
            \zeta_2(t)\Phi_*(\rho)\leq\Phi_*(\rho t)\leq \zeta_3(t)\Phi_*(\rho),~~ \rho, t> 0,
       $$
       $$
            \zeta_2(\|u\|_{\Phi_{*}})\leq\int_\Omega\Phi_{*}(u)dx\leq \zeta_3(\|u\|_{\Phi_*}),~ u\in L_{\Phi_*}(\Omega).
        $$
    \end{prop}



\end{document}